\documentclass[11pt]{article}
\usepackage{latexsym}
\usepackage{amsmath, amsthm}
\usepackage[all,dvips,ps]{xy}
\usepackage{amsfonts}
\usepackage{mathrsfs}

\newtheorem{lem}{Lemma}[section]
\newtheorem{thm}{Theorem}[section]

\newtheorem{ass}{Assumption}[section]
\newtheorem{exa}{Example}[section]
\newtheorem{rem}{Remark}[section]

\newcommand{\tra}{\mathrm{tr}\,}

\newcommand{\Rs}{\mathbb{R}}
\newcommand{\Wc}{\mathcal{W}}
\newcommand{\Uc}{\mathcal{U}}

\newcommand{\Lc}{\mathscr{L}}
\newcommand{\Fs}{\mathscr{F}}
\newcommand{\Rc}{\mathscr{R}}
\newcommand{\Nc}{\mathscr{N}}

\newcommand{\Sn}{{\cal S}^n }
\newcommand{\sn}{{\cal S} }

\newcommand{\bz}{{\bf 0} }

\newcommand{\trace}{{\rm tr\,}}

\newcommand{\bpr}{{\bf Proof.} \hspace{1 em}}
\newcommand{\epr}{ \\ \hspace*{4.5in} $\Box$ }
\newcommand{\beq}{ \begin{equation} }
\newcommand{\eeq}{ \end{equation} }
\newcommand{\bt}{ \begin{tabular} }
\newcommand{\et}{ \end{tabular} }

\begin{document}

\bibliographystyle{plain}
\title{On Farkas Lemma and Dimensional Rigidity of bar Frameworks
 \thanks{Research supported by the Natural Sciences and Engineering
         Research Council of Canada.} }
\vspace{0.3in}
        \author{ A. Y. Alfakih
  \thanks{E-mail: alfakih@uwindsor.ca}
  \\
          Department of Mathematics and Statistics \\
          University of Windsor \\
          Windsor, Ontario N9B 3P4 \\
          Canada
}

\date{\today}
\maketitle

\noindent {\bf AMS classification:} 90C22, 90C25, 52C25, 05C62.

\noindent {\bf Keywords:} Farkas Lemma,  Bar frameworks, dimensional rigidity, 
universal rigidity, facial reduction, semidefinite programming, stress matrices.
\vspace{0.1in}

\begin{abstract}
We present a new semidefinite Farkas lemma involving a side constraint on the rank.
This lemma is then used to present a new proof of a recent characterization, by
Connelly and Gortler \cite{cg14}, of dimensional rigidity of bar frameworks.
\end{abstract}

\section{Introduction}

The celebrated Farkas lemma is at the core of optimization theory. It underpins duality
theory of linear programming, and its semidefinite version plays a key role in strong
duality results of semidefinite programming. As an example of {\em theorems of the alternative},
Farkas lemma establishes the infeasibility of a given linear matrix inequality by exhibiting 
a solution for another linear matrix inequality. In this paper, we present a new semidefinite
Farkas lemma (Theorem \ref{thmmain} below) involving a side constraint on the rank. 
This Farkas lemma is then used to 
provide a new proof of a recent characterization, by Connelly and Gortler \cite{cg14}, of
dimensional rigidity of bar frameworks.   

A {\em bar framework} in $\Rs^r$, denoted by $(G,p)$,
is a simple connected undirected graph $G=(V, E)$ whose nodes are points
$p^1,\ldots,p^n$ in $\Rs^r$; and whose edges are line segments, each joining a pair
of these points. We say that $(G,p)$ is $r$-dimensional if the points 
$p^1,\ldots,p^n$ affinely span $\Rs^r$.

Let $(G,p)$ and $(G,p')$ be two $r$-dimensional and $s$-dimensional
frameworks in $\Rs^r$ and $\Rs^s$ respectively. Then
$(G,p')$ is {\em equivalent} to $(G,p)$ if:
\beq \label{defd}
||{p'}^i-{p'}^j||^2 = ||p^i - p^j ||^2 \quad \text{for each $\{i,j\} \in E(G)$},
\eeq
where $||.||$ denotes the Euclidean norm and $E(G)$ denotes the edge set of $G$.
Moreover, $(G,p')$ is said to be {\em affinely equivalent} to $(G,p)$ if
$(G,p')$ is equivalent to $(G,p)$ and ${p'}^i = A p^i + b$ for all $i=1,\ldots,n$,
where $A$ is an $r \times r$ matrix and $b$ is a vector in $\Rs^r$.
Finally, two $r$-dimensional frameworks $(G,p)$ and $(G,p')$ in $\Rs^r$
are {\em congruent} if:
\beq \label{defd}
||{p'}^i-{p'}^j||^2 = ||p^i - p^j ||^2 \quad \text{for all $i,j=1,\ldots,n$}.
\eeq

An $r$-dimensional framework $(G,p)$ is said to be {\em dimensionally rigid} if no
$s$-dimensional framework $(G,p')$, for any $s \geq r+1$, is equivalent to $(G,p)$.
On the other hand, if every $s$-dimensional framework $(G,p')$, for any $s$, that is
equivalent to $(G,p)$ is in fact congruent to $(G,p)$, then framework $(G,p)$ is said to
be {\em universally rigid}. It turns out that dimensional rigidity and universal rigidity are
closely related.

\begin{thm}[Alfakih \cite{alf07a}]
Let $(G,p)$ be an $r$-dimensional bar framework on $n$ vertices in $\Rs^r$, for some $r
\leq n-2$. Then $(G,p)$ is universally rigid if and only if the following two conditions hold:
\begin{enumerate}
\item $(G,p)$ is dimensionally rigid.
\item There does not exist an $r$-dimensional framework $(G,p')$ in $\Rs^r$
that is affinely equivalent, but not congruent, to $(G,p)$. 
\end{enumerate}
\end{thm}

The notion of a stress matrix plays a key role in the study of universal and dimensional rigidities.
An {\em equilibrium stress} (or simply a stress) of $(G,p)$ is a real-valued function
$\omega$ on $E(G)$ such that:
\beq \label{defw}
\sum_{j:\{i,j\} \in E(G)} \omega_{ij} (p^i - p^j) = \bz \mbox{ for all } i=1,\ldots,n.
\eeq

Let $E(\overline{G})$ denote the edge set of graph $\overline{G}$, the complement graph of $G$. i.e., 
\[
E(\overline{G})= \{ \{i,j\}: i \neq j , \{i,j\} \not \in E(G) \},
\]
and let $\omega =(\omega_{ij})$ be a stress of $(G,p)$. Then the
$n \times n$ symmetric matrix $\Omega$ where
\beq \label{defO}
\Omega_{ij} = \left\{ \begin{array}{ll} -\omega_{ij} & \mbox{if } \{i,j\} \in E(G), \\
                        0   & \mbox{if }  \{i,j\}  \in E(\overline{G}), \\
                   {\displaystyle \sum_{k:\{i,k\} \in E(G)} \omega_{ik}} & \mbox{if } i=j,
                   \end{array} \right.
\eeq
is called the {\em stress matrix}
associated with $\omega$, or a stress matrix
of $(G,p)$. 

The following result provides a sufficient condition for the
dimensional rigidity of a given framework.

\begin{thm}[Alfakih \cite{alf07a}]
\label{thmcon1}
Let $(G,p)$ be an $r$-dimensional bar framework on $n$ vertices in $\Rs^r$, for some $r
\leq n-2$. Then $(G,p)$ is dimensionally rigid if it admits  
a positive semidefinite stress matrix $\Omega$ of rank $n-r-1$.
\end{thm}

\begin{figure}[t]
\setlength{\unitlength}{0.8cm}
\begin{picture}(4,7.5)(-8,-4)
\put(0,0){\circle*{0.2}}
\put(0,2){\circle*{0.2}}
\put( 1,1){\circle*{0.2}}
\put(2, 0){\circle*{0.2}}
\put(-3, -3){\circle*{0.2}}

\put(0,0){\line(0,1){2}}
\put(0,0){\line(1,1){1}}
\put(-3,-3){\line(1,1){3}}
\put(-3,-3){\line(3,5){3}}
\put(-3,-3){\line(5,3){5}}
\put(0,2){\line(1,-1){1}}
\put(1,1){\line(1,-1){1}}
\qbezier(0,2)(1.3,1.3)(2,0)

\put(0.1,-0.4){$1$}
\put(-0.2,2.3){$2$}
\put(1.3,1.1){$3$}
\put(2.2,-0.3){$4$}
\put(-3.4,-2.8){$5$}

\end{picture}
\caption{A $2$-dimensional universally rigid bar framework in the plane, 
    where the set of missing edges is $E(\overline{G}) = \{ \{1,4\}, \{3,5\} \}$.
   It admits a positive semidefinite stress matrix of rank 1 but not of rank 2.
   The edge $\{2,4\}$ is drawn as an arc to make edges $\{2,3\}$
   and $\{3,4\}$ visible.  }
\label{fg1}
\end{figure}
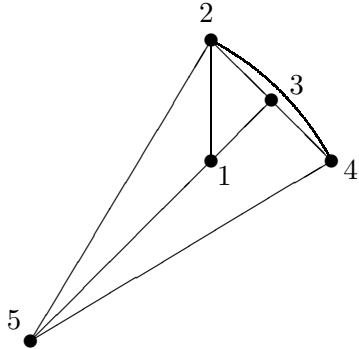

Unfortunately, the sufficient condition in Theorem \ref{thmcon1} is not necessary as was shown by Example 3.1
in \cite{alf07a} (see also Figure \ref{fg1}). 
Recently, Connelly and Gortler \cite{cg14} bridged the gap between necessary and sufficient conditions
for dimensional rigidity. 
Theorem \ref{thmmain2} below is a refined version of their main result in \cite{cg14} concerning
dimensional rigidity.

The remainder of the paper is organized as follows. In Section \ref{secFRFL} we review basic results on
the facial structure of the semidefinite cone and we present our new Farkas lemma. The proof of this lemma  is
based on the Borwein-Wolkowicz facial reduction algorithm \cite{bw81a,bw81b}. In Section \ref{secDRBF}
we review basic results concerning the dimensional rigidity of bar frameworks, and we
use our new Frakas lemma to present a proof of the Connelly-Gortler characterization of dimensional
rigidity in \cite{cg14}. Finally, numerical examples are presented in Section \ref{secNE} to illustrate
the results of the paper.

\subsection{Notation}

For the convenience of the reader, we collect here the notation used throughout the paper.
$I_n$ denotes the identity matrix of order $n$.
$\bz$ denote the zero vector or matrix of appropriate dimension. We denote by $e$
the vector of all 1's in $\Rs^n$, and by $e^i$ we denote the $i$th standard unit vector in $\Rs^n$.
For $i < j$, we let
\beq \label{defFij}
F^{ij}= (e^i-e^j)(e^i-e^j)^T.
\eeq
$||.||$ denotes the Euclidean norm. 
$\Sn$ denotes the space of $n \times n$ symmetric real matrices.
The set of  $n \times n$ symmetric real positive semidefinite (positive definite) matrices is
denoted by $\Sn_+$ ($\Sn_{++}$). We sometimes also use $A \succeq \bz (A \succ \bz)$ to mean that 
$A$ is symmetric positive semidefinite (positive definite).
We denote the relative interior of a set $S$ in $\Sn$ by relint($S$).
For a matrix $A$,
$\Nc (A)$ and $\Rc (A)$ denote, respectively, the null space and the column space (or the range)
of $A$. The trace of $A$ is denoted by $\tra(A)$.
$E(G)$ denotes the edge set of a simple graph $G$, while $E(\overline{G})$ denotes the
edge set of the complement graph of $G$, i.e., 
$E(\overline{G}) = \{ \{i,j\}: i \neq j, \{i,j\} \not \in E(G)\}$.

\section{Facial Reduction and  Farkas Lemma}
\label{secFRFL}

The proof of Theorem \ref{thmmain} below relies on the Borwein-Wolkowicz facial reduction algorithm \cite{bw81a,bw81b}.
Thus we start this section by reviewing definitions and basic facts concerning the facial structure of
the positive semidefinite cone $\Sn_+$. For other applications of facial reduction see
\cite{kri10,kw10,che13}.

\subsection{Facial Structure of $\Sn_+$}
 
A subset $K \in \Sn$, the space of $n \times n$ symmetric real matrices,
 is a {\em cone} if for each $X \in K$ and each $\lambda \geq 0$ we have
$\lambda X \in K$. Let $K$ be a convex cone in $\Sn$. 
A subset $F \subseteq K$ is a {\em face} of $K$ if 
for every $X, Y \in K$ such that $(X+Y) \in F$, it follows that $X \in F$ and $Y \in F$.
A face $F$ of convex cone $K$ is said to be {\em exposed} if there exists an $A \in \Sn$ such
that $F = \{ X \in K: \trace (A X) = 0\}$. A convex cone $K$ is {\em facially exposed} if every
face $F$ of $K$ is exposed.
Let $S$ be a subset of a convex cone $K$, then the intersection of all faces of $K$ containing $S$ is
called the {\em minimal face} of $S$, denoted by face($S$). It is easy to show that
face($S$) is indeed a face of $K$. Moreover, if $S=\{A\}$, we write face($A$) instead
of face$(\{A\})$.

It easy to see that $\Sn_+$, the set of $n \times n$  symmetric positive semidefinite real matrices,
is a closed convex cone. The faces of $\Sn_+$ are well known to be in a one-to-one
correspondence with the subspaces of $\Rs^n$ \cite{bar73,bc75,wsv00}. In fact,
$F$ is a face of $\Sn_+$ if and only if
\beq
F=\{ X \in \Sn_+: \Lc \subseteq \Nc(X) \},
\eeq
for some subspace $\Lc$ of $\Rs^n$, where $\Nc(X)$ denotes the null space of $X$. Moreover,
\beq
\mbox{ relint(} F)=\{ X \in \Sn_+: \Lc = \Nc(X) \}.
\eeq

Thus, faces of $\Sn_+$ are uniquely characterized by their
relative interior. Hence, we have the following theorem. 

\begin{thm}[\cite{bar73,bc75,wsv00}] \label{thmface}
Let $A \in \Sn_+$ of rank $r$ and let 
$A=[W \;\; U] \left[ \begin{array}{cc} \Lambda & 0 \\ 0 & 0 \end{array} \right] 
\left[ \begin{array}{c} W^T \\ U^T \end{array} \right]=W \Lambda W^T$ be the spectral decomposition
of $A$, where $\Lambda$ is the $r \times r$ diagonal matrix consisting of the positive
eigenvalues of $A$. Then 
\begin{eqnarray}
 \mbox{face}(A) & = & \{ X \in \Sn_+: X U = \bz \}, \\
    & = & \{ X \in \Sn_+: X = W Y W^T \mbox{ for some } Y \in \sn^r_+ \}.
\end{eqnarray} 
\end{thm}
Note that $A$ belongs to relint(face($A$)). Let $\Rc(A)$ denote the column space of $A$. Then
for any $B$ in $\Sn_+$ such that 
$\Rc(B) \subset \Rc(A)$, it follows that face($B) \subset $ face($A$). Hence,
if $X \in $ face($A$), then rank $X \leq $ rank $A$. Moreover,
\[
\mbox{if } \Rc(B) = \Rc(A), \mbox{ then face}(B) =  \mbox{face}(A).
\]

\begin{rem}
Observe that face($A$) in Theorem \ref{thmface} has dimension $r(r+1)/2$. More precisely,
face($A$) is isomorphic to $\sn^r_+$. Thus the faces of $\Sn_+$ are isomorphic to
smaller dimensional positive semidefinite cones. Furthermore, it is easy to see that
$\Sn_+$ = face($I_n$).
\end{rem}

The following lemma can be used to provide a characterization of dimensional rigidity.

\begin{lem}\label{lemmain1}
Let $A^1,\ldots,A^m$ be given $n \times n$ symmetric matrices and let $b$ be a given nonzero vector in
$\Rs^m$. Further, let $\Fs=\{ X \in \Sn_+: \tra(X A^i) = b_i$ for $i=1,\ldots,m\}$.
Assume that $X^* \in \Fs$ such that rank $X^*=r$. Then 
there does not exist an $X \in \Fs$ such that rank $X \geq r+1$ 
if and only if $\Fs \subseteq $ face($X^*$).
\end{lem}

\bpr
Assume that $\Fs \subseteq $ face ($X^*$). Then rank $X \leq $ rank ($X^*) = r$ for
all $X \in \Fs$ since $\Rc(X) \subseteq \Rc(X^*)$.

To prove the other direction assume that
rank $(X^*)=r=\max\{$rank $X : X \in \Fs\}$. Let
$X'$ be any matrix in $ \Fs$ and  
let $X = \alpha X^* + (1-\alpha) X'$ for some $\alpha: 0 < \alpha < 1$. Then 
$X \in \Fs$ since $\Fs$ is convex. Furthermore, $\Nc(X)$ = $\Nc(X^*) \cap \Nc(X')$.
Thus, $\Nc(X) \subseteq \Nc(X')$ and $\Nc(X) \subseteq \Nc(X^*)$. 
Hence, rank $X \geq r$. But, $X \in \Fs$. Thus, rank $X=r$. Consequently, $\Nc(X) = \Nc(X^*)$. Therefore,
$\Nc(X^*) \subseteq \Nc(X')$. Hence, $X' \subseteq $ face($X^*$) and thus
$\Fs \subset $ face($X^*$).
\epr 

\begin{rem}\label{rem1}
In fact, it follows from Lemma \ref{lemmain1} that
there does not exist an $X \in \Fs$ such that rank $X \geq r+1$
if and only if face($\Fs) = $ face($X^*$).
This follows since $X^* \in \Fs$ implies that face($X^*) \subseteq $ face($\Fs)$. 
On the other hand, $\Fs \subseteq $ face($X^*$) implies that face($\Fs) \subseteq $ face($X^*$).
\end{rem}

The following lemma plays a key role in this paper.

 \begin{lem} \label{lemBO}
Let $A^1,\ldots, A^m$ be given $n \times n$ symmetric matrices and let $b=(b_i)$ be a given
nonzero vector in $\Rs^m$. Let 
\[
\Fs=\{ X \in \Sn_+: \tra(X A^i) = b_i\mbox{ for } i=1,\ldots,m\}.
\]
Further, let $X^* \in \Fs$ and let $\Uc_j$ be a matrix with full column rank. 
If the following two conditions hold:
\begin{enumerate}
\item $\Fs \subset $ face($\Uc_j \Uc_j^T)$, 
\item There exists $\Omega^j= \sum_{i=1}^m x^j_i A^i$ such that
 $\Uc_j^T \Omega^j \; \Uc_j \succeq \bz, \neq \bz$ and $\tra(\Omega^j X^*) \leq 0$. 
\end{enumerate}
Then 
\beq
\Fs \subset \mbox{face}(\Uc_{j+1} \Uc_{j+1}^T) \subset \mbox{face}(\Uc_{j} \Uc_{j}^T),
\eeq
where $\Wc_j$ is a full column rank matrix such that $\Rc(\Wc_j) = \Nc(\Uc_j \Omega^j \Uc_j)$
and $\Uc_{j+1} = \Uc_j \Wc_j$.
\end{lem}

\bpr 
$\Fs \subset $ face($\Uc_j \Uc_j^T$) implies that 
$\Fs=\{X= \Uc_j Y  \Uc_j^T: Y \succeq \bz,  \tra(X A^i ) = b_i$ for $i=1,\ldots,m\}$.
Then for every $X \in \Fs$ we  have
\[
\tra (X \Omega^j) = \sum_{i=1}^m x^j_{i} \; \tra (\Uc_j Y \Uc_j^T A^i)  
 = \sum_{i=1}^m x^j_{i} b_i
= \sum_{i=1}^m x^j_{i} \; \tra (X^* A^i) 
 = \tra (\Omega^j X^*) \leq \bz.
 \]
But $\tra(X \Omega^j) = \tra (\Uc_j^T \Omega^j \, \Uc_j Y)$. Therefore,
 $ \Uc_j^T \Omega^j \, \Uc_j Y = \bz$ since both $Y \succeq 0$ and $\Uc_j^T \Omega^j \; \Uc_j \succeq \bz$.
 Hence, $Y = \Wc_j Y_j \Wc_j^T$ for some $Y_i \succeq \bz$. Hence, 
 $\Fs=\{X= \Uc_{j+1} Y_j  \Uc_{j+1}^T: Y_j \succeq \bz,  \tra(X A^i) = b_i$ for $i=1,\ldots,m\}$;
 i.e., $\Fs \subset $ face($\Uc_{j+1} \Uc_{j+1}$).
 The result follows since $\Rc(\Uc_{j+1}) \subset \Rc(\Uc_j)$.
\epr

\begin{rem}
In Lemma \ref{lemBO}, suppose that $\Uc_j$ is $n \times s$. Thus,  
face($\Uc_j \Uc_j^T$) is isomorphic to $\sn^s_+$.
Now if rank $(\Uc_j^T \Omega^j \Uc_j) = \delta$, then
$\Wc_j$ is $s \times (s-\delta)$ and hence, 
$\Uc_{j+1}$ is $n \times (s-\delta)$ with full column rank. Consequently,  
face($\Uc_{j+1}\Uc_{j+1}^T$) is isomorphic to $\sn^{(s-\delta)}_+$.
Therefore, the higher the rank of ($\Uc_j^T \Omega^j \Uc_j)$) is, the
larger  the difference between the dimension of 
face($\Uc_{j+1} \Uc_{j+1}^T$) and the dimension of face($\Uc_j \Uc_j^T$) will be.
\end{rem}

\subsection{A New Farkas Lemma}

The following semidefinite Farkas lemma is well known. It is used to establish strong
duality for semidefinite programming under Slater condition (see e.g \cite{lov03}). It will also
be used repeatedly in our proofs. 

\begin{lem} \label{nhomofarkas}
Let $A^1,\ldots, A^m$ be given $n \times n$ symmetric matrices and let $b=(b_i)$ be a given
nonzero vector in $\Rs^m$.  Further, let 
\[
\Fs=\{ X \in \Sn_+: \tra(X A^i) = b_i \mbox{ for }i=1,\ldots,m\}.
\]
Assume that there exists an $X^* \in \Fs$. Then 
exactly one of the following two statement holds:
\begin{enumerate}
\item There exists an $X \in \Fs$ such that $X \succ \bz$,
\item There exists $\Omega = x_1A^1+\cdots+ x_m A^m$ such that $\Omega \succeq \bz, \neq \bz$ and
$\tra(\Omega X^*) \leq 0$.
\end{enumerate}
 \end{lem}

Now are ready to state and prove our new semidefinite Farkas lemma.

\begin{thm} \label{thmmain}
Let $A^1,\ldots, A^m$ be given $n \times n$ symmetric matrices and let $b=(b_i)$ be a given
nonzero vector in $\Rs^m$. Let 
\[
\Fs=\{ X \in \Sn_+: \tra(X A^i) = b_i \mbox{ for }i=1,\ldots,m\}.
\]
Let $\Uc_0$ be a given $n \times s$ matrix with full column rank, and assume that
$\Fs \subset $ face($\Uc_0 \Uc_0^T$).
Let $X^*=\Uc_0 Y^* \Uc_0^T$ be a matrix  in $\Fs$ such that rank $X^* = r$, $ r \leq s-1$. Then 
exactly one of the following two statement holds. 
\begin{enumerate}
\item There exists an $X$ in $\Fs$ such that rank $X \geq r+1$.
\item There exist nonzero matrices $\Omega^0, \Omega^1,\ldots,\Omega^k$, for some $k \leq s-r$, such that: 
\begin{enumerate}
\item $\Omega^j = \sum_{i=1}^m x_i^j A^i$ ($j=0,1,\ldots,k$) for some scalars $x^j_i$,
\item $\Uc_j^T \Omega^j \Uc_j \succeq \bz$ for $j=0,1,\ldots,k$,
\item $\tra(X^* \Omega^j) \leq 0$ for $j=0,1,\ldots,k$,
\item rank $(\Uc_0^T \Omega^0 \Uc_0)$ + rank $(\Uc_1^T \Omega^1 \Uc_1) + \cdots $+  
                                         rank $(\Uc_k^T \Omega^k \Uc_k) = s-r$,
\end{enumerate}
\end{enumerate}
where  $\Uc_1,\ldots, \Uc_{k+1}$, and $\Wc_0,\Wc_1, \ldots,\Wc_{k}$ are full column rank matrices defined as
follows: For $i=0,1,\ldots,k$, 
$\Rc(\Wc_i)= \Nc(\Uc_i^T \Omega^i \Uc_i)$ and 
$\Uc_{i+1}=\Uc_{i} \Wc_i$.
\end{thm}

Before presenting the proof of Theorem \ref{thmmain}, we outline the key idea and intuition behind it.
By Remark  \ref{rem1}, Statement 1 of Theorem \ref{thmmain} does not hold if and only if
face($\Fs )= $ face($X^*$).
Borwein and Wolkowicz \cite{bw81a,bw81b} presented a facial reduction algorithm 
for finding face($\Fs$). At each step of this  algorithm, a smaller dimensional face of $\Sn_+$
containing face($\Fs$) is found. Thus, this algorithm will find matrices 
$\Uc_1, \ldots, \Uc_{k+1}$ such that

\[
\mbox{face}(\Fs) = \mbox{face}(\Uc_{k+1}\Uc_{k+1}^T) \subset \cdots \subset \mbox{face}(\Uc_1\Uc_1^T)
                                                        \subset \mbox{face}(\Uc_0 \Uc_0^T),
\]
where 
$\Rc(\Uc_{k+1}) \subset \cdots \subset \Rc(\Uc_1) \subset \Rc(\Uc_0)$. 
Hence, Statement 1 in Theorem \ref{thmmain} does not hold   
if and only if face($\Uc_{k+1} \Uc_{k+1}^T$) = face($X^*$) if and only if $\Rc(\Uc_{k+1})=\Rc(X^*)$.
   
\bpr
First, we prove that if Statement 1 does not hold, then Statement 2 holds. Therefore,
assume that there does not exist an $X \in \Fs$ such that rank $X \geq r+1$, i.e.,
assume that face($\Fs$)= face($X^*$). Then, 
there does not exist an $s \times s$ matrix $Y \succ \bz$ such that 
$\tra(Y \Uc_0^T A^i \Uc_0) = b_i$ for $i=1,\ldots,m$.
Thus by Lemma \ref{nhomofarkas},
there exists $\Omega^0= \sum_{i=1}^m x_i^0 A^i$ such that $\Uc_0^T \Omega^0 \Uc_0 \succeq \bz, \neq \bz$
and $\tra(X^*\Omega^0) \leq 0$.
If rank ($\Uc_0^T \Omega^0 \Uc_0) = s-r$, then we are done
and $k=0$ in the theorem.
Therefore assume that rank ($\Uc_0^T \Omega^0 \Uc_0) = s-r - \delta_1$, where $\delta_1 \geq 1$,
and let $\Wc_0$ be a full column rank matrix such that $\Rc(\Wc_0)=\Nc(\Uc_0 \Omega^0 \Uc_0)$.
Since $\Fs \subset $ face($\Uc_0 \Uc_0^T$), it follows from Lemma~\ref{lemBO} that 
\[
\Fs \subset \mbox{ face}(\Uc_1 \Uc_1^T) \subset \mbox{ face}(\Uc_0 \Uc_0^T), 
\]
where $\Uc_1=\Uc_0 \Wc_0$ is $n \times (r+\delta_1)$ with full column rank.
Moreover, since face($\Fs$)= face($X^*) \neq $ face($\Uc_1 \Uc_1^T$), 
there does exist $Y_1 \succ \bz$ such that
$\tra(\Uc_1^T A^{i} \Uc_1 Y_1) = b_i$ for all $i=1,\ldots,m$.
Thus, by Lemma \ref{nhomofarkas}, there exists 
$\Omega^1= \sum_{i=1}^m x_i^1 A^i$ such that $\Uc_1^T \Omega^1 \Uc_1 \succeq \bz, \neq \bz$ and
$\tra(X^*\Omega^1) \leq 0$. 
If rank $\Uc_1^T \Omega^1 \Uc_1$ = $\delta_1$, then
rank ($\Uc_0^T \Omega^0 \Uc_0$) + rank ($\Uc_1^T \Omega^1 \Uc_1$) = $s-r$ and we are done and
$k=1$ in the Theorem.  Therefore,
assume that rank ($\Uc_1^T \Omega^1 \Uc_1) = \delta_1 - \delta_2$, where $\delta_1 -1 \geq \delta_2 \geq 1$,
and let $\Wc_1$ be a full column rank matrix such that $\Rc(\Wc_1)=\Nc(\Uc_1 \Omega^1 \Uc_1)$.
Since $\Fs \subset $ face($\Uc_1 \Uc_1^T$), it follows from Lemma \ref{lemBO} that 
\[
\Fs \subset \mbox{ face}(\Uc_2 \Uc_2^T)  \subset \mbox{ face}(\Uc_1 \Uc_1^T)  \subset \mbox{ face}(\Uc_0 \Uc_0^T), 
\]
where $\Uc_2=\Uc_1 \Wc_1$ is $n \times (r+\delta_2)$ with full column rank.

Observe that at each step, a lower dimensional face containing $\Fs$ is obtained. Thus after at most
$s-r$  steps, we must arrive at the case where rank $(\Uc_k \Omega^k \Uc_k) = \delta_k$
and hence Statement 2 holds.

Second, we prove that if Statement 2 holds, then Statement 1 does not hold. Therefore,
assume that $k=0$ in the theorem, i.e., there exists 
$\Omega^0= \sum_{i=1}^m x_i^0 A^i$ such that $\Uc_0^T \Omega^0 \Uc_0 \succeq \bz$, 
rank ($\Uc_0^T \Omega^0 \Uc_0) = s-r$ 
and $\tra(X^*\Omega^0) \leq 0$. 
Since $\Fs \subset $ face($\Uc_0 \Uc_0^T$), it follows from Lemma \ref{lemBO} that 
\[
\Fs \subset \mbox{ face}(\Uc_1 \Uc_1^T), 
\]
where $\Uc_1=\Uc_0 \Wc_0$ is $n \times r$ with full column rank. Hence, 
face($X^*$) = face $(\Uc_1 \Uc_1^T)$ and thus Statement1 does not hold.

Now assume that $k=1$ in the theorem, i.e., there exist 
$\Omega^0 = \sum_{i=1}^m x_i^0 A^i$, $\tra(\Omega^0 X^*) \leq 0$ and 
$\Omega^1 = \sum_{i=1}^m x_i^1 A^i$,  $\tra(\Omega^1 X^*) \leq 0$ 
such that $\Uc_0^T \Omega^0 \Uc_0 \succeq \bz, \neq \bz$ and $\Uc_1^T \Omega^1 \Uc_1 \succeq \bz, \neq \bz$ 
where
rank $\Omega^0=n-r-\delta_1$  and  rank ($\Uc_1^T \Omega^1 \Uc_1$) = $\delta_1$.
Let $\Wc_0$ and $\Wc_1$ be full column rank matrices such that $\Rc(\Wc_0)=\Nc(\Uc_0 \Omega^0 \Uc_0)$
and $\Rc(\Wc_1)=\Nc(\Uc_1 \Omega^1 \Uc_1)$.
Then it follows from Lemma \ref{lemBO} that $\Fs \subset $ face($\Uc_1 \Uc_1^T$) where 
$\Uc_1=\Uc_0 \Wc_0$ is $n \times (r +\delta_1)$ with full column rank.  
Applying Lemma \ref{lemBO} again we have that $\Fs \subset $ face($\Uc_2 \Uc_2^T$) where 
$\Uc_2=\Uc_1 \Wc_1$ is $n \times r$ with full column rank.  
Hence, face($X^*$) = face $(\Uc_2 \Uc_2^T)$ and thus Statement1 does not hold.

Since $\Rc(\Uc_k) \subset \Rc(\Uc_{k-1}) \subset \cdots \subset \Rc(\Uc_1)$, after at most $s-r$ steps
we must have rank ($\Uc_k^T \Omega^k \Uc_k)$ = $\delta_k$. Thus $\Uc_{k+1}=\Uc_k \Wc_k$ is $n \times r$.
Thus face($X^*$) =  face($\Uc_{k+1} \Uc^T_{k+1}$), and thus 
Statement 1 does not hold.
\epr

In Theorem \ref{thmmain}, the assumption that $\Fs \subset $ face($\Uc_0 \Uc_0$) was made in order to make
the application of Theorem \ref{thmmain} to the dimensional rigidity problem straightforward; i.e.,
this assumption was made for the purposes of the paper. Dropping this assumption is equivalent to
setting $\Uc_0=I_n$ since $\Sn_+$ = face($I_n$).  
The following lemma is a restatement of Theorem \ref{thmmain} without
the aforementioned assumption. 

\begin{lem} \label{thmmain2}
Let $A^1,\ldots, A^m$ be given $n \times n$ symmetric matrices and let $b=(b_i)$ be a given
nonzero vector in $\Rs^m$. 
Further, let 
\[
\Fs=\{ X \in \Sn_+: \tra(X A^i) = b_i \mbox{ for }i=1,\ldots,m\}.
\]
Let $X^*$ be a matrix  in $\Fs$ such that rank $X^* = r$. Then 
exactly one of the following two statement holds: 
\begin{enumerate}
\item There exists an $X$ in $\Fs$ such that rank $X \geq r+1$.
\item There exist nonzero matrices $\Omega^0, \Omega^1,\ldots,\Omega^k$, for some $k \leq n-r$, such that: 
\begin{enumerate}
\item $\Omega^j = \sum_{i=1}^m x_i^j A^i$ ($j=0,1,\ldots,k$) for some scalars $x^j_i$,
\item $ \Omega^0 \succeq \bz, \; \Uc_1^T \Omega^1 \Uc_1 \succeq \bz,\; \ldots , \; 
                                                 \Uc_k^T \Omega^k \Uc_k \succeq \bz$,
\item $\tra(X^* \Omega^j) \leq 0$ for $j=0,1,\ldots,k$,
\item rank $\Omega^0 $ + rank $(\Uc_1^T \Omega^1 \Uc_1) + \cdots $+  
                                         rank $(\Uc_k^T \Omega^k \Uc_k) = n-r$,
\end{enumerate}
\end{enumerate}
where  $\Uc_1,\ldots, \Uc_{k+1}$, and $\Wc_0,\Wc_1, \ldots,\Wc_{k}$ are full column rank matrices defined
as follows: For $i=0,1,\ldots,k$, 
$\Rc(\Wc_i)= \Nc(\Uc_i^T \Omega^i \Uc_i)$, and 
$\Uc_{i+1}=\Uc_{i} \Wc_i$ with $\Uc_0=I_n$.
\end{lem}

\section{Dimensional Rigidity of Bar Frameworks}
\label{secDRBF}

In the section, we use Theorem \ref{thmmain} to prove a recent result, by Connelly and Gortler \cite{cg14},
concerning the dimensional rigidity of bar frameworks (or frameworks for short). 
To make the dimensional rigidity problem amenable to semidefinite programming methodology,
we use Gram matrices to represent the configuration of a framework. We start by characterizing 
the set of all frameworks that are equivalent to a given framework $(G,p)$. 
 
\subsection{The Set of Equivalent Frameworks}
\label{secd}

Let $(G,p)$ be an $r$-dimensional framework on $n$ vertices in $\Rs^r$.
The $n \times r$ matrix 
\beq \label{defP}
P= \left[ \begin{array}{c} (p^1)^T \\ \vdots \\ (p^n)^T \end{array} \right]
\eeq
is called the {\em configuration matrix } of $(G,p)$.
We will find it convenient to make the following assumption in the sequel.
Recall that $e$ denotes the vector of all 1's in $\Rs^n$.
\begin{ass} \label{ass1}
$P^Te=0$ for any configuration matrix $P$, i.e., the origin coincides with the centroid
of the points $p^1,\ldots,p^n$.
\end{ass}
In terms of the configuration matrix $P$, the Gram matrix of $(G,p)$ is given by $PP^T$. 
Note that rank $(PP^T) = r$ since $(G,p)$ is $r$-dimensional, i.e., $P$ has full column rank.
Beside being positive semidefinite, Gram matrices of frameworks
are invariant under orthogonal transformation. Moreover, by Assumption \ref{ass1},
Gram matrices of frameworks are also invariant
under translations. Hence, congruent frameworks have the same Gram matrix.
Thus, Gram matrices can be used to characterize all frameworks $(G,p')$ that 
are equivalent to $(G,p)$.

Let $B'$ be the Gram matrix of framework $(G,p')$. Recall the definition of 
matrix $F^{ij}$ in (\ref{defFij}). Then 
\beq 
\tra (F^{ij} B') = B'_{ii} + B'_{jj} - 2 B'_{ij} = || {p'}^i - {p'}^j||^2. 
\eeq
Thus, 
$(G,p')$ is equivalent to $(G,p)$ if and only if
\beq \label{Oprimal}
\tra (F^{ij} B') = ||p^i - p^j ||^2 \mbox{ for all } \{i,j\} \in E(G),
\eeq

The following theorem characterizes the set of all frameworks that are equivalent to
$(G,p)$.

\begin{thm}\label{thmset1}
Let $(G,p)$ be a given $r$-dimensional framework on $n$ nodes in $\Rs^{r}$,
$r \leq n-2$, and let
\beq \label{defF1}
\Fs = \{B' \in \Sn_+: \; B'e=0, \; \tra(F^{ij}B') = ||p^i-p^j||^2  \; \forall \{i, j\} \in E(G)\}.
\eeq
Then $(G,p')$ is an $r'$-dimensional framework that is equivalent to $(G,p)$ if and only
if the Gram matrix of $(G,p')$ belongs to $\Fs$, where $r'$ = rank $B'$.
\end{thm}

The following Theorem is an immediate corollary of Theorem \ref{thmset1} and Lemma~\ref{lemmain1}. 

\begin{thm}\label{drface}
Let $(G,p)$ be a given $r$-dimensional framework on $n$ nodes in $\Rs^{r}$,
$r \leq n-2$, and let
\beq \label{defF11}
\Fs = \{B' \in \Sn_+: \; B'e=0, \; \tra(F^{ij}B') = ||p^i-p^j||^2  \; \forall \{i, j\} \in E(G)\}.
\eeq
Then $(G,p)$ is dimensionally rigid if and only $\Fs \subset $ face($PP^T$), 
where $P$ is the configuration matrix of $(G,p)$.
\end{thm}

Let $V$ be an $n \times (n-1)$ matrix
such that
\beq \label{defV}
V^T e = \bz \;\; \mbox{ and } V^T V = I_{n-1}.
\eeq 
Then, by Theorem \ref{thmface}, it follows that 
$\Fs$ in (\ref{defF1}) is a subset of face($VV^T$).
Thus Theorem \ref{thmset1} can be equivalently stated as follows.
\begin{thm}\label{thmset2}
Let $(G,p)$ be a given $r$-dimensional framework on $n$ nodes in $\Rs^{r}$,
$r \leq n-2$, and let
\beq \label{defF2}
\Fs = \{B'=VYV^T: Y \in \sn^{n-1}_+, \; \tra(Y V^T F^{ij} V ) = ||p^i-p^j||^2  \; \forall \{i, j\} \in E(G)\}.
\eeq
Then $(G,p')$ is an $r'$-dimensional framework that is equivalent to $(G,p)$ if and only
if the Gram matrix of $(G,p')$ belongs to $\Fs$, where $r'$ = rank $Y$.
\end{thm}

\subsection{Quasi-Stress Matrices}

We saw earlier that stress matrices play an important role in the problem of dimensional rigidity.
However, for the purposes of this paper, it will be convenient to introduce
the notion of a quasi-stress matrix.

It is clear from the definition of a stress matrix in (\ref{defO}) that 
the columns of the matrix $[P \;\; e]$ belong to the null space of any stress matrix of $(G,p)$,
where $P$ is the configuration matrix of $(G,p)$.
Hence, the rank of a stress matrix of an $r$-dimensional framework on $n$ vertices is $\leq n-r-1$.
 
An $n \times n$ symmetric matrix $\Omega$ is said to be a {\em quasi-stress matrix} of $(G,p)$ if it satisfies
the following properties: 
\begin{align}
&(a) \;\;\; P^T \Omega P = \bz, \label{quasi2} \\
&(b) \;\;\; \Omega_{ij}=0 \mbox{ for all $\{i,j\} \in E(\overline{G})$}. \nonumber  \\
&(c) \;\;\; \Omega e = \bz. \nonumber
\end{align}

It immediately follows that if a quasi-stress matrix $\Omega$ is positive semidefinite, then $\Omega$ is a stress
matrix since in this case $P^T \Omega P=\bz$ implies that $\Omega P=\bz$. 
For later use we remark here that for any $A \in \Sn$, $Ae=\bz$ if and only if
$A = \sum_{i<j} \omega_{ij} F^{ij}$ for some $\omega_{ij}$'s, where $F^{ij}$ is as defined in (\ref{defFij}).
As a result, any quasi-stress matrix $\Omega$ can be written as
\beq
\Omega = \sum_{\{i,j\} \in E(G)} \omega_{ij} F^{ij} \mbox{ for some scalars } \omega_{ij}. 
\eeq

\subsection{Characterizing Dimensional Rigidity}

A characterization of dimensional rigidity in terms of the minimal face of $PP^T$ was
given in Theorem \ref{drface}. 
Another characterization  can be obtained from Theorem \ref{thmmain}. 
The following theorem is a refined version of  Connelly and Gortler
main result (Corollary 2 in \cite{cg14}) concerning dimensional rigidity.

\begin{thm} \label{thmmain2}
Let $(G,p)$ be an $r$-dimensional framework on $n$ vertices in $\Rs^r$, $r \leq n-2$. Then
$(G,p)$ is dimensionally rigid if and only if there exist nonzero quasi-stress matrices:
$\Omega^0, \Omega^1, \ldots, \Omega^k$, for $k \leq n-r-1$,  such that:  
\begin{enumerate}
\item $\Omega^0 \succeq \bz, \; \Uc_1^T \Omega^1 \Uc_1 \succeq \bz,\; \ldots , \; \Uc_k^T \Omega^k \Uc_k \succeq \bz$,
\item rank $\Omega^0$ + rank $(\Uc_1^T \Omega^1 \Uc_1) + \cdots $+  rank $(\Uc_k^T \Omega^k \Uc_k) = n-r-1$,
\item $P^T \Omega^1 \rho_1 = \bz ,\; \ldots , \; P^T \Omega^k \rho_k = \bz$,
\end{enumerate}
where $\rho_1$, $\Uc_1,\ldots,\Uc_k$ and $\xi_1,\ldots,\xi_k$ are full column rank matrices
defined as follows: $\Rc(\rho_1)= \Nc( \left[ \begin{array}{c} \Omega^0 \\ P^T \\ e^T \end{array} \right])$,
 $\Rc(\xi_i)= \Nc(\rho_i^T \Omega^i \rho_i)$,
$\Uc_i=[P \;\; \rho_i]$ for $i=1,\ldots,k$
and 
$\rho_{i+1}=\rho_i \xi_i$ for all $i=1,\ldots,k-1$.
\end{thm}

\bpr
framework $(G,p)$ is dimensionally rigid if and only if $\Fs \subset $ face($PP^T)$, i.e.,
if and only if there does not exist a $B' \in \Fs$ such that rank $B' \geq r+1$, where
$\Fs$ is defined in (\ref{defF2}). Note that $\Fs \subset $ face($VV^T$), where
$V$ is as defined in (\ref{defV}). Therefore,
it follows from  Theroem \ref{thmmain} that $(G,p)$ is dimensionally rigid
if and only if there exist nonzero  matrices $\Omega^0,\ldots, \Omega^k$, for some $k \leq n-1-r$,
 such that: 
\begin{enumerate}
\item $\Omega^l = \sum_{\{i,j\} \in E(G)} \omega^l_{ij} F^{ij}$ ( $l=0,1,\ldots,k$ ) for some scalars $\omega_{ij}^l$.
\item $V^T \Omega^0 V \succeq \bz, \;\Uc_l^T \Omega^l \Uc_l \succeq \bz$ for $l=1,\ldots,k$,
\item $\tra(PP^T \Omega^l) \leq 0$ for $l=0,1,\ldots,k$,
\item rank $(V^T \Omega^0 V)$ + rank $(\Uc_1^T \Omega^1 \Uc_1) + \cdots $+  
                                         rank $(\Uc_k^T \Omega^k \Uc_k) = n-1-r$,
\end{enumerate}
where  $\Uc_1,\ldots, \Uc_{k+1}$, and $\Wc_0,\Wc_1, \ldots,\Wc_{k}$ are full column rank matrices such that
for $i=0,1,\ldots,k$, we  have
$\Rc(\Wc_i)= \Nc(\Uc_i^T \Omega^i \Uc_i)$, and 
$\Uc_{i+1}=\Uc_{i} \Wc_i$, $\Uc_0=V$.

Now it follows from the definition of $V$ in (\ref{defV}) that
\[
VV^T=I_n - \frac{ee^T}{n}.
\] 
Let $\overline{\Omega}^0=V^T \Omega^0 V$, then $\Omega^0=V \overline{\Omega}^0 V^T$ since $\Omega^0e=\bz$.
Therefore, $\overline{\Omega}^0 \succeq \bz$ if and only if $\Omega^0 \succeq \bz$.
Thus, $\Omega^0$ is a stress matrix of $(G,p)$ and hence $\Omega^0 P =\bz$.
Moreover, rank $\Omega^0$ = rank $(V^T \Omega^0 V)$. 

Also, it is easy to see that if $x \in \Nc(\Omega^0)$, then $V^Tx \in \Nc(V^T \Omega^0 V)$. Thus,
by the definition of $\rho_1$, we have that $\Rc([P \;\; e \;\; \rho_1]) = \Nc(\Omega^0)$.
Thus $\Rc([V^TP \;\;\; V^T\rho_1]) = \Nc(V^T \Omega^0V)$ = $\Wc_0$. Hence,
$\Uc_1=V \Wc_0=[P \;\; \rho_1]$ since $P^Te=\bz$ and $\rho_1 e = \bz$.

On the other hand, since
\[
\Uc_1^T \Omega^1 \Uc_1 = \left[ \begin{array}{cc} P^T \Omega^1 P & P^T \Omega^1 \rho_1 \\
        \rho_1^T \Omega^1 P & \rho_1^T \Omega^1 \rho_1 \end{array} \right] \succeq \bz, 
\]
and since $\tra(P^T \Omega^1 P) \leq 0$, it follows that
\[
P^T \Omega^1 P = \bz \mbox{ and hence } P^T \Omega^1 \xi = \bz.
\]
Therefore, $\Omega^1$ is a quasi-stress matrix of $(G,p)$. Moreover, since
\[
\Rc(\Wc_1) = \Nc( \Uc_1^T \Omega^1 \Uc_1 ) = \Nc(\left[ \begin{array}{cc} \bz &  \bz \\
        \bz & \rho_1^T \Omega^1 \rho_1 \end{array} \right],
\]
it follows that $\Wc_1= \left[ \begin{array}{cc} I_r & \bz \\ \bz & \xi_1 \end{array} \right]$.
Hence, $\Uc_2=\Uc_1 \Wc_1=[P \;\; \rho_1 \xi_1]= [P \;\; \rho_2]$.
The rest of the proof for $\Omega^2,\ldots,\Omega^k$ proceeds in an analogous fashion to the proof for
$\Omega^1$.
\epr

We end this section with the following observation regarding the computation of
$\Omega^1,\ldots,\Omega^k$. While the matrices $\{F^{ij}: \{i,j\} \in E(G)\}$ are
linearly independent, in the second and subsequent steps of the Borwein-Wolkowicz 
facial reduction algorithm, the matrices  $\{\Uc_l^T F^{ij} \Uc_l: \{i,j\} \in E(G)\}$
may become linearly dependent. Thus some of the distance constraints in the definition of
$\Fs$, namely the constraints 
\[
\tra(Y_l \Uc_l^T F^{ij} \Uc_l) = ||p^i-p^j||^2: \{i,j\} \in E(G)
\]
may become redundant. This can be used to our advantage. Since if the constraint corresponding to edge 
$\{\hat{i},\hat{j}\}$ is redundant, then $\Omega^l_{\hat{i},\hat{j}}$ is 0.
For instance, in Example \ref{ex1}, the framework $(G,p)$ (see Figure \ref{fg1}) has a clique
$\{p^2,p^3,p^4\}$ where $p^2,p^3$ and $p^4$ are collinear. As a result,
\[
\Uc_1^T F^{24} \Uc_1 = 4 \, \Uc_1^T F^{23} \Uc_1 = 4 \, \Uc_1^T F^{34} \Uc_1.
\] 
Hence, we may choose $\Omega^1_{24}=\Omega^1_{34}=0$.

\section{Numerical Examples}
\label{secNE}

To illustrate the Borwein-Wolkowicz facial reduction algorithm used in the proof of Theorem \ref{thmmain},
and consequently in Theorem \ref{thmmain2}, we present two numerical examples.

\begin{exa} \label{excg}
Consider the bar framework $(G,p)$ in Figure \ref{ladder} given in \cite{cg14}. Its
configuration matrix and a corresponding  positive semidefinite stress matrix $\Omega^0$ are given by
\[
P=\left[ \begin{array}{rr} -1 & -2 \\ -1 & 2 \\ 1 & 2 \\ 1 & -2 \\ -1 & 0 \\ 1 & 0 \end{array} \right] 
\mbox{ and }
\Omega^0 =\left[ \begin{array}{rr} 1 & 0 \\ 1 & 0 \\ 0 & 1 \\ 0 & 1 \\ -2 & 0 \\ 0 & -2 \end{array} \right]
\left[ \begin{array}{rrrrrr} 1 & 1 & 0 & 0 & -2 & 0 \\ 0 & 0 & 1 & 1& 0 & -2 \end{array} \right].
\]
Thus $\delta_1=n-r-1-$ rank $\Omega^0 = 1$. 
Then $\rho_1$, the matrix whose columns form a basis of
null space $\left[ \begin{array}{c} \Omega^0 \\ P^T \\ e^T \end{array} \right]$,
is given by
\[
\rho_1= \left[ \begin{array}{r} -1 \\ 1 \\ -1 \\ 1 \\ 0 \\ 0 \end{array} \right].
\mbox{ Hence, }
\Uc_1= [P \;\; \rho_1] =\left[ \begin{array}{rrr} -1 & -2 & -1\\ -1 & 2 & 1 \\ 1 & 2 & -1 \\ 
                                1 & -2 & 1 \\ -1 & 0 & 0 \\ 1 & 0 & 0 \end{array} \right].
\]
Thus 
\[
\Fs = \{B' :
 B' = \Uc_1 Y_1 \Uc_1^T: Y_1 \in \sn^{3}_+, \tra (\Uc_1^T F^{ij} \Uc_1 Y_1) = ||p^i-p^j||^2 \; \forall \; \{i,j\} \in E(G) \},
\]
where,
$\Uc_1^T F^{12} \Uc_1 = 4 \; \Uc_1^T F^{15} \Uc_1 =4 \; \Uc_1^T F^{25} \Uc_1 = 
4 \left[ \begin{array}{rrr} 0 & 0 & 0 \\  0 & 4 & 2 \\ 0 & 2 & 1 \end{array} \right]$,
$\Uc_1^T F^{34} \Uc_1  = 4 \; \Uc_1^T F^{36} \Uc_1 =  4 \;  \Uc_1^T F^{46} \Uc_1=
4 \left[ \begin{array}{rrr} 0 & 0 & 0 \\  0 & 4 & -2 \\ 0 & -2 & 1 \end{array} \right]$,
$\Uc_1^T F^{14} \Uc_1  = 
\left[ \begin{array}{rrr} 4 & 0 & 4 \\  0 & 0 & 0 \\ 4 & 0 & 4 \end{array} \right]$,
$\Uc_1^T F^{23} \Uc_1  = 
\left[ \begin{array}{rrr} 4 & 0 & -4 \\  0 & 0 & 0 \\ -4 & 0 & 4 \end{array} \right]$,
and $\Uc_1^T F^{56} \Uc_1  = 
\left[ \begin{array}{rrr} 4 & 0 & 0 \\  0 & 0 & 0 \\ 0 & 0 & 0 \end{array} \right]$.

\begin{figure}[t]
\setlength{\unitlength}{0.8cm}
\begin{picture}(4,7.5)(-8,-4)
\put(-1,-2){\circle*{0.2}}
\put(-1,2){\circle*{0.2}}
\put( 1,2){\circle*{0.2}}
\put(1, -2){\circle*{0.2}}
\put(-1, 0){\circle*{0.2}}
\put(1, 0){\circle*{0.2}}

\put(-1,-2){\line(1,0){2}}
\put(-1,-2){\line(0,1){2}}
\put(-1,0){\line(1,0){2}}
\put(-1,0){\line(0,1){2}}
\put(-1,2){\line(1,0){2}}
\put(1,-2){\line(0,1){2}}
\put(1,0){\line(0,1){2}}
\qbezier(-1,-2)(-1.35,0)(-1,2)
\qbezier(1,-2)(1.35,0)(1,2)

\put(-1.4,-2.4){$1$}
\put(-1.4,2.3){$2$}
\put(1.2,2.3){$3$}
\put(1.2,-2.3){$4$}
\put(-1.5,-0.2){$5$}
\put(1.3,-0.2){$6$}

\end{picture}
\caption{The dimensionally but not universally rigid bar framework of Example \ref{excg}. 
   It admits a positive semidefinite stress matrix of rank 2 but not of rank 3.
   The edges $\{1,2\}$ and $\{3,4\}$ are drawn as arcs to make edges $\{1,5\}$,
    $\{2,5\}$, $\{3,6\}$ and $\{4,6\}$ visible.  }
\label{ladder}
\end{figure}
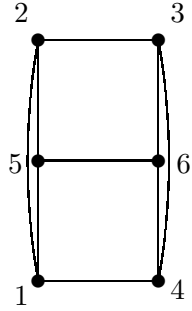

Note that since $\Uc_1^T F^{12} \Uc_1 = 4 \; \Uc_1^T F^{15} \Uc_1 =4 \; \Uc_1^T F^{25} \Uc_1$
and since $\Uc_1^T F^{34} \Uc_1 = 4 \; \Uc_1^T F^{36} \Uc_1 =4 \; \Uc_1^T F^{46} \Uc_1$,
we only include the distance constraints for edges $\{1,2\}$ and $\{3,4\}$. The distance constraints
corresponding to edges $\{1,5\}$, $\{2,5\}$, $\{3,6\}$ and $\{4,6\}$ are redundant.
Thus $Y_1=(y_{ij})$ must satisfy

\[
\begin{array}{ll}
4 \; y_{22} + 4 \; y_{23} + y_{33} & = 4 \\
4 \; y_{22} - 4 \; y_{23} + y_{33} & = 4 \\
y_{11} - 2 \; y_{13} + y_{33} & = 1 \\
y_{11} + 2 \; y_{13} + y_{33} & = 1 \\
4 \; y_{11} & = 4.
\end{array}
\] 
Hence, $y_{11}=y_{22}=1$, $y_{13} = y_{23} = y_{33} = 0$, and $y_{12}= \alpha$ is a free parameter.
Thus $Y_1$ is a function of $\alpha$, and it is given by 
\[
Y_1 (\alpha) = \left[ \begin{array}{rrr} 1 & \alpha & 0 \\  \alpha & 1 & 0 \\ 0 & 0 & 0 \end{array} \right].
\]
Thus, $Y_1(\alpha) $ is positive semidefinite iff 
\[
-1 \leq \alpha \leq 1.
\]
Furthermore, rank $Y_1(\alpha) \leq 2$. Thus $\Fs \; \cap $ relint(face $(\Uc_1\Uc_1^T)) = \emptyset$;
i.e., there does not exist $Y_1 \succeq \bz$ such that $\tra(Y_1 \Uc_1 f^{ij} \Uc_1) = ||p^i-p^j||^2$ for
all $\{i,j\} \in E(G)$.

Now let $\Omega^1=(\omega_{ij})$. 
Since the distance constraints
corresponding to edges $\{1,5\}$, $\{2,5\}$, $\{3,6\}$ and $\{4,6\}$ are redundant, we set
$\omega_{15}=\omega_{25}=\omega_{36}=\omega_{46}=0$.
Then 
$P^T \Omega^1 P = \bz$ and $P^T \Omega^1 \rho_1 = \bz$ imply that 
\begin{eqnarray*}
\omega_{14} + \omega_{23} + \omega_{56} & = & 0, \\
\omega_{14} - \omega_{23} & = & 0, \\
\omega_{12} - \omega_{34} & = & 0, \\
\omega_{12} + \omega_{34} & = & 0.
\end{eqnarray*}
Thus, $\omega_{12}=\omega_{34}=0$ and 
$\omega_{56} = -2 \omega_{14} = -2 \omega_{23}$.
Also, $\rho_1 \Omega^1 \rho_1$ is nonzero positive semidefinite if
$\omega_{12}+\omega_{34}+\omega_{14}+\omega_{23} > 0$.
Therefore, set $\omega_{14} = \omega_{23} = 1$ and $\omega_{56}=-2$. Hence,
\[
\Omega^1 =\left[ \begin{array}{rrrrrr} 1 & 0  & 0 & -1 & 0 & 0 \\  0 & 1 & -1 & 0 & 0 & 0 \\ 
                                    0 & -1  & 1 & 0  &  0 & 0 \\ -1& 0 & 0  & 1 & 0 & 0 \\
                                    0 & 0   & 0 & 0  &  -2 & 2 \\ 0& 0 & 0  & 0 & 2 & -2 \end{array} \right].
\]
Then 
\[ \Uc_1^T \Omega^1 \Uc_1 = \left[ \begin{array}{rrc}  0 & 0 & 0 \\ 0 & 0 & 0\\
                                                  0 & 0 & \rho_1^T \Omega_1 \rho_1 \end{array} \right]
=  \left[ \begin{array}{rrr}  0  & 0 & 0 \\ 0 & 0 & 0 \\ 0 & 0 & 8  \end{array} \right] \succeq \bz.
\]
Hence,  $\rho_1 \Omega^1 \rho_1 = 8$ is nonsingular. Moreover,
$\Rc( \left[ \begin{array}{cc} I_2 \\ 0  \end{array} \right])= \Nc(\Uc_1^T \Omega^1 \Uc_1)$.
Hence, $\Uc_2 = \Uc_1 \left[ \begin{array}{cc} I_2 \\ 0  \end{array} \right]$=
 $[P \;\; \rho_1] \left[ \begin{array}{cc} I_2 \\ 0  \end{array} \right]=P$.
Therefore, $\Fs \subseteq $ face($PP^T$) and
rank $\Omega^0$ + rank ($\Uc_1^T \Omega^1 \Uc_1$) = 2+1 = 3 = $n-r-1$.
\end{exa}

\begin{exa} \label{ex1}
Consider the bar framework $(G,p)$ in Figure \ref{fg1}. Its
configuration matrix and a corresponding positive semidefinite stress matrix $\Omega^0$ are given by
\[
P=\left[ \begin{array}{rr} 0 & 0 \\ 0 & 2 \\ 1 & 1 \\ 2 & 0 \\ -3 & -3  \end{array} \right] 
\mbox{ and }
\Omega^0 =\left[ \begin{array}{rr} 0  \\ 1  \\ -2 \\  1 \\ 0  \end{array} \right]
\left[ \begin{array}{rrrrr} 0 & 1  & -2 & 1 & 0  \end{array} \right].
\]
Thus $\delta_1=n-r-1-$ rank $\Omega^0 = 1$. 
Then $\rho_1$, the matrix whose columns form a basis of
null space $\left[ \begin{array}{c} \Omega^0 \\ P^T \\ e^T \end{array} \right]$,
is given by
\[
\rho_1= \left[ \begin{array}{r} -4 \\ 1 \\ 1 \\ 1 \\ 1  \end{array} \right].
\mbox{ Hence, }
\Uc_1= [P \;\; \rho_1] =\left[ \begin{array}{rrr} 0 & 0 & -4\\ 0 & 2 & 1 \\ 1 & 1 & 1 \\ 
                                2 & 0 & 1 \\ -3 & -3 & 1 \end{array} \right].
\]
Thus 
\[
\Fs = \{B' :
 B' = \Uc_1 Y_1 \Uc_1^T: Y_1 \in \sn^{3}_+, \tra (\Uc_1^T F^{ij} \Uc_1 Y_1) = 
 ||p^i-p^j||^2 \; \forall \; \{i,j\} \in E(G) \},
\]
where,
$\Uc_1^T F^{24} \Uc_1 = 4 \; \Uc_1^T F^{23} \Uc_1 =4 \; \Uc_1^T F^{34} \Uc_1 = 
4 \left[ \begin{array}{rrr} 1 & -1 & 0 \\  -1 & 1 & 0 \\ 0 & 0 & 0 \end{array} \right]$,
$\Uc_1^T F^{12} \Uc_1  = 
\left[ \begin{array}{rrr} 0 & 0 & 0 \\  0 & 4 & 10 \\ 0 & 10 & 25 \end{array} \right]$,
$\Uc_1^T F^{13} \Uc_1  = 
\left[ \begin{array}{rrr} 1 & 1 & 5 \\  1 & 1 & 5 \\ 1 & 5 & 25 \end{array} \right]$,
$\Uc_1^T F^{15} \Uc_1  = 
\left[ \begin{array}{rrr} 9 & 9 & -15 \\  9 & 9 & -15 \\ -15 & -15 & 25 \end{array} \right]$,
$\Uc_1^T F^{25} \Uc_1  = 
\left[ \begin{array}{rrr} 9 & 15 & 0 \\  15 & 25 & 0 \\ 0 & 0 & 0 \end{array} \right]$.
and $\Uc_1^T F^{45} \Uc_1  = 
\left[ \begin{array}{rrr} 25 & 15 & 0 \\  15 & 9 & 0 \\ 0 & 0 & 0 \end{array} \right]$.

Note that since $\Uc_1^T F^{24} \Uc_1 = 4 \; \Uc_1^T F^{23} \Uc_1 =4 \; \Uc_1^T F^{34} \Uc_1$, 
we only include the distance constraint for edge $\{2,3\}$. The distance constraints
corresponding to edges $\{2,4\}$ and $\{3,4\}$ are redundant.
Thus $Y_1=(y_{ij})$ must satisfy

\[
\begin{array}{ll}
y_{11} - 2 \; y_{12} + y_{22} & = 2 \\
25 \; y_{11} + 30 \; y_{12} + 9 \; y_{22} & = 34 \\
9 \; y_{11} + 30 \; y_{12} + 25 \; y_{22} & = 34 
\end{array}
\]
and
\[
\begin{array}{ll}
9 \; y_{11} + 18 \; y_{12} + 9 \; y_{22} -30 \; y_{13} - 30 \; y_{23} + 25 \; y_{33} & = 18 \\
y_{11} + 2 \; y_{12} +  y_{22} +10 \; y_{13} +10 \; y_{23} + 25 \; y_{33} & = 2 \\
4 \; y_{22} +20 \; y_{23}  + 25 \; y_{33} & = 4
\end{array}
\] 
Hence, $Y_1$ is unique (recall that $(G,p)$ is universally rigid), and it is given by 
\[
Y_1 = \left[ \begin{array}{rrr} 1 & 0 & 0 \\  0 & 1 & 0 \\ 0 & 0 & 0 \end{array} \right].
\]

Furthermore, rank $Y_1 = 2$. Thus $\Fs \; \cap $ relint(face $(\Uc_1\Uc_1^T)) = \emptyset$. 

Now let $\Omega^1=(\omega_{ij})$. 
Since the distance constraints
corresponding to edges $\{2,4\}$ and $\{3,4\}$ are redundant, we set $\omega_{24}=\omega_{34}=0$.
Then 
$P^T \Omega^1 P = \bz$ and $P^T \Omega^1 \rho_1 = \bz$ imply that 
\begin{eqnarray*}
 \omega_{13} - 3 \;  \omega_{15}  & = & 0, \\
 2 \; \omega_{12} + \omega_{13} - 3 \; \omega_{15} & = & 0, \\
 \omega_{13} + 9 \; \omega_{15} + \omega_{23} + 9 \; \omega_{25} + 25 \; \omega_{45}  & = & 0, \\
 \omega_{13} + 9 \; \omega_{15} - \omega_{23} + 15 \; \omega_{25} + 15 \; \omega_{45}  & = & 0, \\
4 \; \omega_{12} + \omega_{13} + 9 \;\omega_{15} + \omega_{23} + 25 \; \omega_{25} + 9 \; \omega_{45}  & = & 0. \\
\end{eqnarray*}
Moreover, we require  $\rho_1^T \Omega^1 \rho_1 = 25(\omega_{12}+\omega_{13}+\omega_{15})$ to be positive semidefinite. 
 Hence, $\omega_{12}=0$,  $\omega_{13}=24$,  $\omega_{15}=8$,
 $\omega_{23}=6$,  and $\omega_{25}= \omega_{45}=-3$; i.e., 
\[
\Omega^1 =\left[ \begin{array}{rrrrr} 32 & 0  & -24 & 0 & -8 \\  0 & 3 & -6 & 0 & 3 \\ 
                                    -24 & -6  & 30 & 0  &  0 \\ 0 & 0  & 0 & -3 & 3 \\
                                    -8   & 3 & 0  &  3 & 2 \end{array} \right].
\]
Then 
\[ \Uc_1^T \Omega^1 \Uc_1 = \left[ \begin{array}{rrc}  0 & 0 & 0 \\ 0 & 0 & 0\\
                                                  0 & 0 & \rho_1^T \Omega_1 \rho_1 \end{array} \right]
=  \left[ \begin{array}{rrr}  0  & 0 & 0 \\ 0 & 0 & 0 \\ 0 & 0 & 800  \end{array} \right] \succeq \bz.
\]
Hence,  $\rho_1 \Omega^1 \rho_1 = 800$ is nonsingular. Moreover,
$\Rc( \left[ \begin{array}{cc} I_2 \\ 0  \end{array} \right])= \Nc(\Uc_1^T \Omega^1 \Uc_1)$.
Hence, $\Uc_2 = \Uc_1 \left[ \begin{array}{cc} I_2 \\ 0  \end{array} \right]$=
 $[P \;\; \rho_1] \left[ \begin{array}{cc} I_2 \\ 0  \end{array} \right]=P$.
Therefore, $\Fs \subseteq $ face($PP^T$) and
rank $\Omega^0$ + rank ($\Uc_1^T \Omega^1 \Uc_1$) = 2+1 = 3 = $n-r-1$.
\end{exa}



\end{document}